\documentclass[reqno]{amsart}
\usepackage{amssymb,setspace}
\usepackage{ifpdf}
\ifpdf
  \usepackage[hyperindex,pagebackref]{hyperref}
\else
  \expandafter\ifx\csname dvipdfm\endcsname\relax
    \usepackage[hypertex,hyperindex,pagebackref]{hyperref}
  \else
    \usepackage[dvipdfm,hyperindex,pagebackref]{hyperref}
  \fi
\fi
\theoremstyle{plain}
\newtheorem{theorem}{Theorem}

\newtheorem{proposition}{Proposition}
\theoremstyle{definition}
\newtheorem{definition}{Definition}
\DeclareMathOperator{\td}{d\mspace{-2mu}}
\allowdisplaybreaks[4]
\numberwithin{equation}{section}

\begin{document}

\title[The $q$-Genocchi numbers and polynomials and applications]
{On the $\boldsymbol{q}$-Genocchi numbers and polynomials with weight zero and their applications}

\author[S. Arac\i]{Serkan Arac\i}
\address[Arac\i]{Department of Mathematics, Faculty of Science and Arts, University of Gaziantep, 27310 Gaziantep, TURKEY}
\email{\href{mailto: S. Arac\i <mtsrkn@hotmail.com>}{mtsrkn@hotmail.com}}

\author[M. A\c{c}ikg\"oz]{Mehmet A\c{c}ikg\"oz}
\address[A\c{c}ikg\"oz]{Department of Mathematics, Faculty of Science and Arts, University of Gaziantep, 27310 Gaziantep, TURKEY}
\email{\href{mailto: M. Acikgoz <acikgoz@gantep.edu.tr>}{acikgoz@gantep.edu.tr}}

\author[F. Qi]{Feng Qi}
\address[Qi]{School of Mathematics and Informatics, Henan Polytechnic University, Jiaozuo City, Henan Province, 454010, China; Department of Mathematics, School of Science, Tianjin Polytechnic University, Tianjin City, 300387, China}
\email{\href{mailto: F. Qi <qifeng618@gmail.com>}{qifeng618@gmail.com}, \href{mailto: F. Qi <qifeng618@hotmail.com>}{qifeng618@hotmail.com}, \href{mailto: F. Qi <qifeng618@qq.com>}{qifeng618@qq.com}}
\urladdr{\url{http://qifeng618.wordpress.com}}



\begin{abstract}
In this paper, the authors deal with the $q$-Genocchi numbers and polynomials with weight
zero. They discover some interesting relations via the $p$-adic $q$-integral
on $\mathbb{Z}_{p}$ and familiar basis Bernstein polynomials. Finally, the authors show that the $p$-adic $\log$ gamma functions are associated with the $q$-Genocchi numbers and polynomials with weight zero.
\end{abstract}

\maketitle

\section{Preliminaries}

Let $p$ be an odd prime number. Denote the ring of the $p$-adic rational integers by $\mathbb{Z}_{p}$, the field of rational numbers by $\mathbb{Q}$, the field of the $p$-adic rational numbers by $\mathbb{Q}_{p}$, and the completion of algebraic closure of $\mathbb{Q}_{p}$ by $\mathbb{C}_{p}$, respectively. Let $\mathbb{N}$ be the set of positive integers and $\mathbb{N}^{\ast}=\mathbb{N}\cup\{0\}$ the set of all non-negative integers. The $p$-adic absolute value is defined by
\begin{equation}
|p|_{p}=\frac1p.
\end{equation}
Assume $|q-1|_{p}<1$ is an indeterminate number in the sense that either $q\in\mathbb{C}$ or $q\in\mathbb{C}_p$. A $q$-analogue of $x$ may be defined by
\begin{equation}
[x]_{q}=\frac{1-q^{x}}{1-q}
\end{equation}
satisfying $\lim_{q\to 1}[x]_{q}=x$.
\par
A function $f$ is said to be uniformly differentiable at a point $a\in\mathbb{Z}_{p}$ if the divided difference
\begin{equation*}
F_{f}(x,y)=\frac{f(x)-f(y)}{x-y}
\end{equation*}
converges to $f'(a)$ as $(x,y)\to (a,a)$. The class of all the uniformly differentiable functions is denoted by $UD(\mathbb{Z}_{p})$.
\par
For $f\in UD(\mathbb{Z}_{p})$, the $p$-adic $q$-analogue of Riemann sum for $f$ was defined by
\begin{equation}\label{Reieman-sum-q-qnqlogue}
\frac{1}{[p^{n}]_q}\sum_{0\le\xi <p^{n}}f(\xi)q^{\xi}
=\sum_{0\le\xi <p^{n}}f(\xi)\mu_{q}\bigl(\xi+p^{n}\mathbb{Z}_{p}\bigr)
\end{equation}
in~\cite{Araci-Genocchi-Kim0, Araci-Genocchi-kim7}, where $n\in\mathbb{N}$.
The integral of $f$ on $\mathbb{Z}_{p}$ is defined as the limit of~\eqref{Reieman-sum-q-qnqlogue} as $n$ tends to $\infty$, if it exists, and represented by
\begin{equation}
I_{q}(f)=\int_{\mathbb{Z}_{p}}f(\xi) \td\mu_{q}(\xi). \label{Genocchi-eq1}
\end{equation}
The bosonic integral and the fermionic $p$-adic integral on $\mathbb{Z}_{p}$ are defined respectively by
\begin{equation}
I_{1}(f)=\lim_{q\to 1}I_{q}(f)
\end{equation}
and
\begin{equation}
I_{-q}(f)=\lim_{q\to -q}I_{q}(f). \label{Genocchi-eq2}
\end{equation}
For a prime $p$ and a positive integer $d$ with $(p,d)=1$, set
\begin{gather*}
X=X_{d}=\lim_{\overleftarrow{n}}\mathbb{Z}/dp^{n}\mathbb{Z}, \quad X_{1}=\mathbb{Z}_{p}, \\
X^{\ast}=\bigcup_{\substack{(a,p)=1\\0<a<dp}}a+dp\mathbb{Z}_{p},
\end{gather*}
and
\begin{equation*}
a+dp^{n}\mathbb{Z}_{p}=\bigl\{ x\in X\mid x\equiv a\mod dp^{n}\bigr\},
\end{equation*}
where $a\in\mathbb{Z}$ satisfies $0\le a<dp^{n}$ and $n\in\mathbb{N}$.

\section{Main results}

In~\cite{Araci-Genocchi-1, Araci-Genocchi-2}, Arac\'i, A\c{c}ikg\"oz, and Seo considered the $q$-Genocchi polynomials with weight $\alpha$ in the form
\begin{equation}
\frac{\widetilde{G}_{n+1,q}^{(\alpha)}(x)}{n+1}
=\int_{\mathbb{Z}_{p}}[ x+\xi]_{q^{\alpha}}^{n}\td\mu_{-q}(\xi),
\label{Genocchi-eq3}
\end{equation}
where $\widetilde{G}_{n+1,q}^{(\alpha)}=\widetilde{G}_{n+1,q}^{(\alpha)}(0)$ is called the $q$-Genocchi numbers with weight $\alpha$.
Taking $\alpha=0$ in~\eqref{Genocchi-eq3}, we easily see that
\begin{equation}
\frac{\widetilde{G}_{n+1,q}}{n+1}\triangleq\frac{\widetilde{G}_{n+1,q}^{(0)}}{n+1}=\int_{\mathbb{Z}_{p}}\xi^{n}\td\mu_{-q}(\xi),   \label{Genocchi-eq4}
\end{equation}
where $\widetilde{G}_{n,q}$ are called the $q$-Genocchi numbers and polynomials
with weight $0$. From~\eqref{Genocchi-eq4}, it is simple to see
\begin{equation}
\sum_{n=0}^{\infty}\widetilde{G}_{n,q}\frac{t^{n}}{n!}=t\int_{\mathbb{Z}_{p}}e^{\xi t}\td\mu_{-q}(\xi) .  \label{Genocchi-eq5}
\end{equation}
By~\eqref{Genocchi-eq2}, we have
\begin{equation}
q^{n}I_{-q}(f_{n})+(-1)^{n-1}I_{-q}(f)
=[2]_{q}\sum_{0\le\ell<n}q^{\ell}(-1)^{n-1-\ell}f(\ell),  \label{Genocchi-eq6}
\end{equation}
where $f_{n}(x)=f(x+n)$ and $n\in\mathbb{N}$. See~\cite{Araci-Genocchi-kim16, Araci-Genocchi-kim15, Araci-Genocchi-kim17}. Taking $n=1$ in~\eqref{Genocchi-eq6} leads to the well-known equality
\begin{equation}
qI_{-q}(f_{1})+I_{-q}(f)=[2]_{q}f(0),  \label{Genocchi-eq7}
\end{equation}
When setting $f(x)=e^{xt}$ in~\eqref{Genocchi-eq7}, we find
\begin{equation}
\sum_{n=0}^{\infty}\widetilde{G}_{n,q}\frac{t^{n}}{n!}=\frac{[2]_{q}t}{qe^{t}+1}.  \label{Genocchi-eq8}
\end{equation}
By~\eqref{Genocchi-eq8}, we obtain the $q$-Genocchi polynomials with weight $0$ as follows
\begin{equation}\label{Genocchi-eq9}
\sum_{n=0}^{\infty}\widetilde{G}_{n,q}(x)\frac{t^{n}}{n!}=
\frac{[2]_{q}t}{qe^{t}+1}e^{xt}.
\end{equation}
By~\eqref{Genocchi-eq9}, we see that
\begin{equation*}
\sum_{n\geq 0}\widetilde{G}_{n,q}(x)\frac{t^{n}}{n!}=t\frac{
1-\bigl(-q^{-1}\bigr)}{e^{t}-\bigl(-q^{-1}\bigr)}e^{xt}=t\sum_{n\geq
0}H_{n}\bigl(-q^{-1},x\bigr) \frac{t^{n}}{n!}.
\end{equation*}
By equating coefficients of $t^{n}$ on both sides of the above equality, we derive the following theorem.

\begin{theorem}
For $n\in\mathbb{N}$, we have
\begin{equation*}
\frac{\widetilde{G}_{n+1,q}(x)}{n+1}=H_{n}\bigl(-q^{-1},x\bigr),
\end{equation*}
where $H_{n}\bigl(-q^{-1},x\bigr)$ are the $n$-th Frobenius-Euler polynomials.
\end{theorem}

By~\eqref{Genocchi-eq7}, we discover that
\begin{align*}
[2]_{q}\sum_{n=0}^{\infty}x^{n}\frac{t^{n}}{n!} &=q\int_{\mathbb{Z}_{p}}e^{(x+\xi+1) t}\td\mu_{-q}(\xi)+\int_{\mathbb{Z}_{p}}e^{(x+\xi) t}\td\mu_{-q}(\xi)  \\
&=\sum_{n=0}^{\infty}\biggl[q\int_{\mathbb{Z}_{p}}(x+\xi+1)^{n}\td\mu_{-q}(\xi)
+\int_{\mathbb{Z}_{p}}(x+\xi)^{n}\td\mu_{-q}(\xi)\biggr]\frac{t^{n}}{n!} \\
&=\sum_{n=0}^{\infty}\bigl[qH_{n}\bigl(-q^{-1},x+1\bigr)+H_{n}\bigl(-q^{-1},x\bigr)\bigr]\frac{t^{n}}{n!}.
\end{align*}
Equating coefficients of $\frac{t^{n}}{n!}$ on both sides above equation, we deduce the following theorem.

\begin{theorem}\label{Genocchi-eq10-thm}
For $n\in\mathbb{N}$, the identity
\begin{equation}
qH_{n}\bigl(-q^{-1},x+1\bigr)+H_{n}\bigl(-q^{-1},x\bigr)=[ 2
]_{q}x^{n}  \label{Genocchi-eq10}
\end{equation}
is valid.
\end{theorem}

In particular, when letting $q=1$, the identity~\eqref{Genocchi-eq10} becomes
\begin{equation}
G_{n}(x+1)+G_{n}(x)=2nx^{n-1},
\end{equation}
where $G_{n}(x)$ are called the Genocchi polynomials.
\par
If we substitute $x=0$ into~\eqref{Genocchi-eq10}, then Theorem~\ref{Genocchi-eq10-thm} can be rewritten as Theorem~\ref{Genocchi-eq13-thm} below.

\begin{theorem}\label{Genocchi-eq13-thm}
The identity
\begin{equation}
q{\widetilde{G}_{n,q}(1)}+{\widetilde{G}_{n,q}}=\begin{cases}[2]_{q},& n=1\\
0, &  n\ne 1
\end{cases}  \label{Genocchi-eq13}
\end{equation}
is true, where $\widetilde{G}_{n,q}$ are called the Genocchi numbers and
polynomials with weight $0$.
\end{theorem}

When we substitute $x$ by $1-x$ and $q$ by $q^{-1}$ in~\eqref{Genocchi-eq9}, it follows that
\begin{multline*}
\sum_{n=0}^{\infty}\widetilde{G}_{n,q^{-1}}(1-x)\frac{t^{n}}{n!
}=t\frac{1+q^{-1}}{q^{-1}e^{t}+1}e^{(1-x) t}
=\frac{1+q}{e^{t}+q}e^{t}e^{xt} \\*
=-\frac{[2]_{q}(-t)}{qe^{-t}+1}e^{(-t) x}
=\sum_{n=0}^{\infty}(-1)^{n+1}\widetilde{G}_{n,q}(x)\frac{t^{n}}{n!}.
\end{multline*}
From this, we procure symmetric properties of this type polynomials.

\begin{theorem}
The following identity holds
\begin{equation}
\widetilde{G}_{n,q^{-1}}(1-x)=(-1)^{n+1}\widetilde{
G}_{n,q}(x).  \label{Genocchi-eq11}
\end{equation}
\end{theorem}

By using~\eqref{Genocchi-eq3} for $\alpha=0$ and the binomial theorem, we readily obtain that
\begin{align*}
\frac{\widetilde{G}_{n+1,q}(x)}{n+1}&=\int_{\mathbb{Z}_{p}}(x+\xi)^{n}\td\mu_{-q}(\xi) \\
&=\sum_{k=0}^{n}\binom{n}{k}\biggl[\int_{\mathbb{Z}_{p}}\xi^{k}\td\mu_{-q}(\xi)\biggr] x^{n-k} \\
&=\sum_{k=0}^{n}\binom{n}{k}\frac{\widetilde{G}_{k+1,q}}{k+1}x^{n-k}.
\end{align*}
Further using
$$
\frac{n+1}{k+1}\binom{n}{k}=\binom{n+1}{k+1},
$$
we obtain
$$
\widetilde{G}_{n+1,q}(x)=\sum_{k=0}^n\binom{n+1}{k+1}\widetilde{G}_{k+1,q}x^{n-k} =\sum_{k=1}^{n+1}\binom{n+1}k\widetilde{G}_{k,q}x^{n+1-k}.
$$
Thus, we get the following conclusion.

\begin{theorem}
The identity
\begin{equation}
\widetilde{G}_{n,q}(x)=\sum_{k=0}^{n}\binom{n}{k}\widetilde{G}
_{k,q}x^{n-k} \label{Genocchi-eq12}
\end{equation}
is true, where the usual convention of replacing $\bigl(\widetilde{G}_{q}\bigr)^{n}$ by $\widetilde{G}_{n,q}$ is used.
\end{theorem}

Combining~\eqref{Genocchi-eq13} with~\eqref{Genocchi-eq12} leads to the following proposition.

\begin{proposition}
The identity
\begin{equation}
\widetilde{G}_{0,q}=0\quad \text{and}\quad {q\bigl(\widetilde{G}_{q}+1\bigr)^{n}} +{\widetilde{G}_{n,q}}=\begin{cases}[2]_{q},& n=1\\
0, &  n\ne 1
\end{cases}  \label{Genocchi-eq14}
\end{equation}
is true, where the usual convention of replacing $\bigl(\widetilde{G}_{q}\bigr)^{n}$ by $\widetilde{G}_{n,q}$ is used.
\end{proposition}

From~\eqref{Genocchi-eq12}, it follows that
\begin{align*}
q^{2}\widetilde{G}_{n+1,q}(2) &=q^{2}\bigl(\widetilde{G}_{q}+1+1\bigr)^{n+1}\\
&=q^{2}\sum_{k=0}^{n+1}\binom{n+1}{k}\bigl(\widetilde{G}_{q}+1\bigr)^{k} \\
&=(n+1) q^{2}\bigl(\widetilde{G}_{q}+1\bigr)
^{1}+q\sum_{k=2}^{n+1}\binom{n+1}{k}q\bigl(\widetilde{G}_{q}+1\bigr)^{k}\\
&=(n+1) q\bigl([2]_{q}-\widetilde{G}_{1,q}\bigr)
-q\sum_{k=2}^{n+1}\binom{n+1}{k}\widetilde{G}_{k,q} \\
&=(n+1) q[2]_{q}-\Biggl[q\sum_{k=2}^{n+1}\binom{n+1}{k}\widetilde{G}_{k,q} +(n+1)q\widetilde{G}_{1,q}\Biggr]\\
&=(n+1) q[2]_{q}-q\sum_{k=0}^{n+1}\binom{n+1}{k}\widetilde{G}_{k,q}\\
&=(n+1) q[2]_{q}-q\bigl(\widetilde{G}_{q}+1\bigr)^{n+1} \\
&=(n+1) q[2]_{q}+\widetilde{G}_{n+1,q}
\end{align*}
for $n>1$. Therefore, we deduce the following proposition.

\begin{proposition}
For $n>1$,
\begin{equation}
\widetilde{G}_{n+1,q}(2)=\frac{(n+1)}{q}[2]_{q}+\frac{1}{q^{2}}\widetilde{G}_{n+1,q}.
\label{Genocchi-eq15}
\end{equation}
\end{proposition}

By virtue of~\eqref{Genocchi-eq2}, \eqref{Genocchi-eq11}, and~\eqref{Genocchi-eq15}, we find
\begin{multline*}
(n+1)\int_{\mathbb{Z}_{p}}(1-\xi)^{n}\td\mu_{-q}(\xi)=(n+1) (-1)^{n}\int_{\mathbb{Z}_{p}}(\xi -1)^{n}\td\mu
_{-q}(\xi)  \\
=(-1)^{n}\widetilde{G}_{n+1,q}(-1)=\widetilde{G}_{n+1,q^{-1}}(2)
=(n+1) [2]_{q}+q^{2}\widetilde{G}_{n+1,q^{-1}}.
\end{multline*}
As a result, we may concluded Theorem~\ref{Genocchi-eq16-thm} below.

\begin{theorem}\label{Genocchi-eq16-thm}
The identity
\begin{equation}
\int_{\mathbb{Z}_{p}}(1-\xi)^{n}\td\mu_{-q}(\xi)=[2]
_{q}+q^{2}\frac{\widetilde{G}_{n+1,q^{-1}}}{n+1}  \label{Genocchi-eq16}
\end{equation}
is valid.
\end{theorem}

 Let $UD(\mathbb{Z}_{p})$ be the space of continuous functions on $\mathbb{Z}_{p}$. For $f\in UD(\mathbb{Z}_{p})$, the $p$-adic analogue of Bernstein operator for $f$ is defined by
\begin{equation*}
\boldsymbol{B}_{n}(f,x)=\sum_{k=0}^{n}f\biggl(\frac{k}{n}\biggr) B_{k,n}(x)
=\sum_{k=0}^{n}f\biggl(\frac{k}{n}\biggr)\binom{n}{k}x^{k}(1-x)^{n-k},
\end{equation*}
where $n,k\in\mathbb{N}^{\ast}$ and the $p$-adic Bernstein polynomials of degree $n$ is defined by
\begin{equation}
B_{k,n}(x)=\binom{n}{k}x^{k}(1-x)^{n-k},\quad
x\in\mathbb{Z}_{p}.\label{Genocchi-eq19}
\end{equation}
See~\cite{Araci-Genocchi-3, Araci-Genocchi-Kim-3, Araci-Genocchi-k-2, Araci-Genocchi-Kim-4}.
Via the $p$-adic $q$-integral on $\mathbb{Z}_{p}$ and Bernstein polynomials in~\eqref{Genocchi-eq19}, we can obtain that
\begin{align*}
I_{1} &=\int_{\mathbb{Z}_{p}}B_{k,n}(\xi) \td\mu_{-q}(\xi)\\
&=\binom{n}{k}\int_{\mathbb{Z}_{p}}\xi^{k}(1-\xi)^{n-k}\td\mu_{-q}(\xi) \\
&=\binom{n}{k}\sum_{\ell=0}^{n-k}\binom{n-k}{\ell}(-1)^{\ell} \biggl[\int_{\mathbb{Z}_{p}}\xi^{\ell+k}\td\mu_{-q}(\xi)\biggr] \\
&=\binom{n}{k}\sum_{\ell=0}^{n-k}\binom{n-k}{\ell}(-1)^{\ell}\frac{
\widetilde{G}_{\ell+k+1,q}}{\ell+k+1}.
\end{align*}
On the other hand, by symmetric properties of Bernstein polynomials, we have
\begin{align*}
I_{2} &=\int_{\mathbb{Z}_{p}}B_{n-k,n}(1-\xi) \td\mu_{-q}(\xi)\\
&=\binom{n}{k}\sum_{s=0}^{k}\binom{k}{s}(-1)^{k+s}\int_{\mathbb{Z}_{p}}(1-\xi)^{n+s}\td\mu_{-q}(x)\\
&=\binom{n}{k}\sum_{s=0}^{k}\binom{k}{s}(-1)^{k+s} \biggl([2]_{q}+q^{2}\frac{\widetilde{G}_{n+s+1,q^{-1}}}{n+s+1}\biggr) \\
&=
\begin{cases}\displaystyle
[2]_{q}+q^{2}\frac{\widetilde{G}_{n+s+1,q^{-1}}}{n+s+1},& k=0\\\displaystyle
\binom{n}{k}\sum_{s=0}^{k}\binom{k}{s}(-1)^{k+s}\biggl([2]_{q}+q^{2}\frac{\widetilde{G}
_{n+s+1,q^{-1}}}{n+s+1}\biggr), &  k\ne 0.
\end{cases}
\end{align*}
Equating $I_{1}$ and $I_{2}$ yields Theorem~\ref{thm7-araci-qi} below.

\begin{theorem}\label{thm7-araci-qi}
The following identity holds:
\begin{equation*}
\sum_{\ell=0}^{n-k}\binom{n-k}{\ell}(-1)^{\ell}\frac{\widetilde{G}_{\ell+k+1,q}}{\ell+k+1}=
\begin{cases}\displaystyle
[2]_{q}+q^{2}\frac{\widetilde{G}_{n+s+1,q^{-1}}}{n+s+1},& k=0;\\\displaystyle
\sum_{s=0}^{k}\binom{k}{s}(-1)^{k+s}\biggl([2]_{q}+q^{2} \frac{\widetilde{G}_{n+s+1,q^{-1}}}{n+s+1}\biggr), &  k\ne 0.
\end{cases}
\end{equation*}
\end{theorem}

The $p$-adic $q$-integral on $\mathbb{Z}_{p}$ of the product of several Bernstein polynomials can be calculated as
\begin{align*}
I_{3} &=\int_{\mathbb{Z}_{p}}\prod_{s=1}^{m}B_{k,n_{s}}(\xi) \td\mu_{-q}(
\xi) \\
&=\prod_{s=1}^{m}\binom{n_{s}}{k}\int_{\mathbb{Z}_{p}} \xi^{mk}(1-\xi)^{n_{1}+\dotsm+n_{m}-mk}\td\mu_{-q}(\xi)\\
&=\prod_{s=1}^{m}\binom{n_{s}}{k}\sum_{\ell=0}^{n_{1}+\dotsm+n_{m}-mk}
\binom{n_{1}+\dotsm+n_{m}-mk}{\ell}(-1)^{\ell} \biggl[\int_{\mathbb{Z}_{p}}\xi^{\ell+mk}\td\mu_{-q}(\xi)\biggr]\\
&=\prod_{s=1}^{m}\binom{n_{s}}{k}\sum_{\ell=0}^{n_{1}+\dotsm+n_{m}-mk}
\binom{n_{1}+\dotsm+n_{m}-mk}{\ell}(-1)^{\ell}\frac{\widetilde{G}
_{\ell+mk+1,q^{-1}}}{\ell+mk+1}.
\end{align*}
On the other hand, by symmetric properties of Bernstein polynomials and~\eqref{Genocchi-eq16}, we have
\begin{align*}
I_{4} &=\int_{\mathbb{Z}_{p}}\prod_{s=1}^{m}B_{n_{s}-k,n_{s}}(1-\xi) \td\mu
_{-q}(\xi) \\
&=\binom{n}{k}\sum_{\ell=0}^{mk}\binom{mk}{\ell}(
-1)^{mk+\ell}\int_{\mathbb{Z}_{p}}(1-\xi)^{n_{1}+\dotsm+n_{m}+\ell}\td\mu_{-q}(\xi) \\
&=\prod_{s=1}^{m}\binom{n_{s}}{k}\sum_{\ell=0}^{mk}\binom{mk}{\ell}(
-1)^{mk+\ell}\biggl([2]_{q}+q^{2}\frac{\widetilde{G}
_{n_{1}+\dotsm+n_{m}+\ell+1,q^{-1}}}{n_{1}+\dotsm+n_{m}+\ell+1}\biggr) \\
&=\begin{cases}\displaystyle
[2]_{q}+q^{2}\frac{\widetilde{G}
_{n_{1}+\dotsm+n_{m}+1,q^{-1}}}{n_{1}+\dotsm+n_{m}+1}, & k=0\\\displaystyle
\prod_{s=1}^{m}\binom{n_{s}}{k}\sum_{\ell=0}^{mk}\binom{mk}{\ell}(
-1)^{mk+\ell}\biggl([2]_{q}+q^{2}\frac{\widetilde{G}
_{n_{1}+\dotsm+n_{m}+\ell+1,q^{-1}}}{n_{1}+\dotsm+n_{m}+\ell+1}\biggr), & k\ne 0.
\end{cases}
\end{align*}
Equating $I_{3}$ and $I_{4}$ results in an interesting identity for $q$-analogue of Genocchi polynomials with weight $0$.

\begin{theorem}
The identity
\begin{multline*}
\sum_{\ell=0}^{n_{1}+\dotsm+n_{m}-mk}\binom{n_{1}+\dotsm+n_{m}-mk}{\ell}(-1)
^{\ell}\frac{\widetilde{G}_{\ell+mk+1,q^{-1}}}{\ell+mk+1}\\
=\begin{cases}\displaystyle
[2]_{q}+q^{2}\frac{\widetilde{G}_{n_{1}+\dotsm+n_{m}+1,q^{-1}}}{n_{1}+\dotsm+n_{m}+1}, & k=0\\\displaystyle
\sum_{\ell=0}^{mk}\binom{mk}{\ell}(-1)^{mk+\ell} \biggl([2]_{q}+q^{2}\frac{\widetilde{G}_{n_{1}+\dotsm+n_{m}+\ell+1,q^{-1}}} {n_{1}+\dotsm+n_{m}+\ell+1}\biggr), &   k\ne 0
\end{cases}
\end{multline*}
is true.
\end{theorem}

\section{Other identities}

In this section, we consider Kim's $p$-adic $q$-$\log$ gamma functions related to the $q$-analogue of Genocchi polynomials.

\begin{definition}[\cite{Araci-Genocchi-Kim-1, Araci-Genocchi-Kim0}]
For $x\in\mathbb{C}_{p}\setminus\mathbb{Z}_{p}$,
\begin{equation*}
(1+x)\log(1+x)=x+\sum_{n=1}^{\infty}\frac{(-1)^{n+1}}{n(n+1)}x^{n+1}.
\end{equation*}
\end{definition}

Kim's $p$-adic locally analytic function on $x\in\mathbb{C}_{p}\setminus\mathbb{Z}_{p}$ can be defined as follows.

\begin{definition}[\cite{Araci-Genocchi-Kim-1, Araci-Genocchi-Kim0}]
For $x\in\mathbb{C}_{p}\setminus\mathbb{Z}_{p}$,
\begin{equation*}
G_{p,q}(x)=\int_{\mathbb{Z}_{p}}[ x+\xi]_{q}(\log[ x+\xi]_{q}-1)
\td\mu_{-q}(\xi).
\end{equation*}
\end{definition}

If $q\to 1$, then
\begin{equation} \label{Genocchi-eq17}
G_{p,1}(x)\triangleq G_{p}(x)=\int_{\mathbb{Z}_{p}}(x+\xi)[\log(x+\xi)-1] \td\mu_{-q}(\xi).
\end{equation}
Replacing $x$ by $\frac{\xi}{x}$ in~\eqref{Genocchi-eq17} leads to
\begin{equation}
(x+\xi)[\log(x+\xi)-1]=(x+\xi)\log x+\sum_{n=1}^{\infty}\frac{(-1)^{n+1}}{
n(n+1)}\frac{\xi^{n+1}}{x^{n}}-x.  \label{Genocchi-eq18}
\end{equation}
From~\eqref{Genocchi-eq17} and~\eqref{Genocchi-eq18}, we can establish an interesting formula~\eqref{araci-qi-final-eq} which is useful for studying in the theory of the $p$-adic analysis and the analytic number.

\begin{theorem}
For $x\in\mathbb{C}_{p}\setminus\mathbb{Z}_{p}$,
\begin{equation}\label{araci-qi-final-eq}
G_{p}(x)=\biggl(x+\frac{\widetilde{G}_{2,q}}{2}\biggr)\log x+\sum_{n=1}^{\infty}\frac{(-1)^{n+1}}{n(n+1)
(n+2)}\frac{\widetilde{G}_{n+2,q}}{x^{n}}-x.
\end{equation}
\end{theorem}

\end{document}